\documentclass[12pt,a4paper]{amsart}

\def\frak{\mathfrak}

\newtheorem{theorem}{Theorem}[section]
\newtheorem{lemma}[theorem]{Lemma}
\newtheorem{corollary}[theorem]{Corollary}
\newtheorem{conjecture}[theorem]{Conjecture}
\newtheorem{proposition}[theorem]{Proposition}

\theoremstyle{definition}
\newtheorem{definition}[theorem]{Definition}

\theoremstyle{remark}
\newtheorem{remark}[theorem]{Remark}

\numberwithin{equation}{section}

\def\br{\mathbb R}
\def\bc{\mathbb C}

\def\bz{\mathbb Z}

\def\mand{\quad\mbox{and}\quad}

\def\p{\partial}
\def\Ad{{\rm Ad}}

\def\bp{\begin{pmatrix}}
\def\ep{\end{pmatrix}}
\def\op{{\rm op}}

\begin{document}

\title[Equivalence of domains III]{Equivalence of domains arising from duality of orbits on flag manifolds III}

\author{Toshihiko MATSUKI}
\address{Department of Mathematics\\
        Faculty of Science\\
        Kyoto University\\
        Kyoto 606-8502, Japan}
\email{matsuki@math.kyoto-u.ac.jp}
\date{}

\begin{abstract}
In \cite{GM1}, we defined a $G_\br$-$K_\bc$ invariant subset $C(S)$ of $G_\bc$ for each $K_\bc$-orbit $S$ on every flag manifold $G_\bc/P$ and conjectured that the connected component $C(S)_0$ of the identity would be equal to the Akhiezer-Gindikin domain $D$ if $S$ is of nonholomorphic type. This conjecture was proved for closed $S$ in \cite{WZ1,WZ2,FH,M8} and for open $S$ in \cite{M8}. It was proved for the other orbits in \cite{M9} when $G_\br$ is of non-Hermitian type. In this paper, we prove the conjecture for an arbitrary non-closed $K_\bc$-orbit when $G_\br$ is of Hermitian type. Thus the conjecture is completely solved affirmatively.
\end{abstract}

\maketitle

\section{Introduction}

Let $G_\bc$ be a connected complex semisimple Lie group and $G_\br$ a connected real form of $G_\bc$. Let $K_\bc$ be the complexification in $G_\bc$ of a maximal compact subgroup $K$ of $G_\br$. Let $X=G_\bc/P$ be a flag manifold of $G_\bc$ where $P$ is an arbitrary parabolic subgroup of $G_\bc$. Then there exists a natural one-to-one correspondence between the set of $K_\bc$-orbits $S$ and the set of $G_\br$-orbits $S'$ on $X$ given by the condition:
\begin{equation}
S\leftrightarrow S'\Longleftrightarrow S\cap S'\mbox{ is non-empty and compact} \tag{1.1}
\end{equation}
(\cite{M4}). For each $K_\bc$-orbit $S$ we defined in \cite{GM1} a subset $C(S)$ of $G_\bc$ by
$$C(S)=\{x\in G_\bc\mid xS\cap S'\mbox{ is non-empty and compact}\}$$
where $S'$ is the $G_\br$-orbit on $X$ given by (1.1).

Akhiezer and Gindikin defined a domain $D/K_\bc$ in $G_\bc/K_\bc$ as follows (\cite{AG}). Let $\frak{g}_\br=\frak{k}\oplus\frak{m}$ denote the Cartan decomposition of $\frak{g}_\br={\rm Lie}(G_\br)$ with respect to $K$. Let $\frak{t}$ be a maximal abelian subspace in $i\frak{m}$. Put
$$\frak{t}^+=\{Y\in\frak{t}\mid |\alpha(Y)|<{\pi\over 2} \mbox{ for all }\alpha\in\Sigma\}$$
where $\Sigma$ is the restricted root system of ${\frak g}_\bc$
with respect to $\frak{t}$. Then $D$ is defined by
$$D=G_\br(\exp\frak{t}^+)K_\bc.$$

We conjectured the following in \cite{GM1}.

\begin{conjecture} {\rm (Conjecture 1.6 in \cite{GM1})} \ Suppose that $X=G_\bc/P$ is not $K_\bc$-homogeneous. Then we will have $C(S)_0=D$ for all $K_\bc$-orbits $S$ of nonholomorphic type on $X$. Here $C(S)_0$ is the connected component of $C(S)$ containing the identity.
\end{conjecture}

\begin{remark} \ When $G_\br$ is of Hermitian type, there exist two special closed $K_\bc$-orbits $S_1=K_\bc B/B=Q/B$ and $S_2=K_\bc w_0B/B=w_0Q/B$ on the full flag manifold $G_\bc/B$ where $Q=K_\bc B$ is the usual maximal parabolic subgroup of $G_\bc$ defined by a nontrivial central element in $i\frak{k}$ and $w_0$ is the longest element in the Weyl group. For each parabolic subgroup $P$ containing the Borel subgroup $B$, two closed $K_\bc$-orbits $S_1P$ and $S_2P$ on $G_\bc/P$ are called of holomorphic type and all the other $K_\bc$-orbits are called of nonholomorphic type. Especially all the non-closed $K_\bc$-orbits are defined to be of nonholomorphic type.

When $G_\br$ is of non-Hermitian type, we define that all the $K_\bc$-orbits are of nonholomorphic type.
\end{remark}

Let $S_\op$ denote the unique open dense $K_\bc$-$B$ double coset in $G_\bc$. Then $S'_\op$ is the unique closed $G_\br$-$B$ double coset in $G_\bc$. In this case we see that
$$C(S_\op)=\{x\in G_\bc\mid xS_\op\supset S'_\op\}.$$
It follows easily that $C(S_\op)$ is a Stein manifold (c.f. \cite{GM1}, \cite{H}). The connected component $C(S_\op)_0$ is often called the Iwasawa domain.

The inclusion
$$D\subset C(S_\op)_0$$
was proved in \cite{H}. (Later \cite{Ms} gave a proof without complex analysis.) On the other hand, it was proved in \cite{GM1} Proposition 8.1 and Proposition 8.3 that
$C(S_\op)_0\subset C(S)_0$
for all $K_\bc$-$P$ double cosets $S$ for any $P$. So we have the inclusion
\begin{equation}
D\subset C(S)_0. \tag{1.2}
\end{equation}

Hence we have only to prove the converse inclusion
\begin{equation}
C(S)_0\subset D \tag{1.3}
\end{equation}
for $K_\bc$-orbits $S$ of nonholomorphic type in Conjecture 1.1.

If $S$ is closed in $G_\bc$, then we can write
$$C(S)=\{x\in G_\bc\mid xS\subset S'\}.$$
So the connected component $C(S)_0$ is essentially equal to the cycle space introduced in \cite{WW}. For Hermitian cases the inclusion (1.3) for closed $S$ was proved in \cite{WZ1} and \cite{WZ2}. For non-Hermitian cases it was proved in \cite{FH} and \cite{M8}.

When $S$ is the open $K_\bc$-$P$ double coset in $G_\bc$, the inclusion (1.3) was proved in \cite{M8} for arbitrary $P$ generalizing the result in \cite{B}.

Recently the inclusion (1.3) was proved in \cite{M9} for an arbitrary orbit $S$ when $G_\br$ is of non-Hermitian type. So the remaining problem was to prove (1.3) for non-closed and non-open orbits when $G_\br$ is of Hermitian type.

In this paper we solve this problem.

In the next section we prove the following theorem.

\begin{theorem} \ Suppose that $G_\br$ is of Hermitian type and let $S$ be a non-closed $K_\bc$-$P$ double coset in $G_\bc$. Then there exist $K_\bc$-$B$ double cosets $\widetilde{S}_1$ and $\widetilde{S}_2$ contained in the boundary $\p S=S^{cl}-S$ of $S$ such that
$$x(\widetilde{S}_1\cup \widetilde{S}_2)^{cl}\cap {S'_0}^{cl}\ne\phi$$
for all the elements $x$ in the boudary of $D$. Here $S_0$ denote the dense $K_\bc$-$B$ double coset in $S$.
\end{theorem}

\begin{remark} \ By computations of examples it seems that $\widetilde{S_1}$ and $\widetilde{S_2}$ are always distinct $K_\bc$-orbits. But we do not need this distinctness.
\end{remark}

\begin{corollary} \ Suppose that $G_\br$ is of Hermitian type and let $S$ be a non-closed $K_\bc$-$P$ double coset in $G_\bc$. Then $C(S)_0=D$.
\end{corollary}

\begin{proof} \ Let $S_0$ be as in Theorem 1.3. Let $\Psi$ denote the set of the simple roots in the  positive root system for $B$. For each $\alpha\in\Psi$ we define a parabolic subgroup
$$P_\alpha=B\sqcup Bw_\alpha B$$
of $G_\bc$. By \cite{GM2} Lemma 2 we can take a sequence $\{\alpha_1,\ldots,\alpha_m\}$ of simple roots such that
$$\dim_\bc S_0P_{\alpha_1}\cdots P_{\alpha_k}=\dim_\bc S_0+k$$
for $k=0,\ldots,m={\rm codim}_\bc S_0$. Then it is shown in \cite{M9} Theorem 1.2 that
\begin{equation}
x\in C(S)\cap D^{cl}\Longrightarrow xS^{cl}\cap S'_\op P_{\alpha_m}\cdots P_{\alpha_1}= xS\cap S'_0. \tag{1.4}
\end{equation}

Let $x$ be an element in the boundary of $D$. Then it follows from Theorem 1.3 that
$$x(\p S)\cap {S'_0}^{cl}\ne\phi.$$
If $x$ is also contained in $C(S)$, then it follows from (1.4) that
$$x(\p S)\cap S'_\op P_{\alpha_m}\cdots P_{\alpha_1}=\phi.$$
Since ${S'_0}^{cl}$ is contained in the closed set $S'_\op P_{\alpha_m}\cdots P_{\alpha_1}$, we have
$$x(\p S)\cap {S'_0}^{cl}=\phi,$$
a contradiction. Hence $x\notin C(S)$. Thus we have proved
$C(S)_0\subset D$.
\end{proof}

Section 3 is devoted to the explicit computation of the case where $G_\br=Sp(2,\br)$. We use Proposition 3.2 in the proof of Lemma 2.4 in Section 2. Another simple example of $SU(2,1)$-case is explicitly computed in \cite{M8} Example 1.5.

\bigskip
\noindent {\bf Acknowledgement.} \ The author would like to express his heartily thanks to S.
Gindikin for valuable suggestions and encouragements.

\section{Proof of Theorem 1.3}

Let $\frak{j}$ be a maximal abelian subspace of $i\frak{k}$. Let $\Delta$ denote the root system of the pair $(\frak{g}_\bc, \frak{j})$. Since $G_\br$ is a group of Hermitian type, there exists a nontrivial central element $Z$ of $i\frak{k}$ and we can write
$$\frak{g}_\bc=\frak{k}_\bc\oplus \frak{n}\oplus \overline{\frak{n}}$$
where $\Delta^+_n=\{\alpha\in\Delta \mid \alpha(Z)>0\}$ and  $\frak{n}=\bigoplus_{\alpha\in\Delta^+_n} \frak{g}_\bc(\frak{j},\alpha)$. Let $Q$ be the maximal parabolic subgroup of $G_\bc$ defined by
$$Q=K_\bc\exp\frak{n}.$$
Let $\Delta^+$ be a positive system of $\Delta$ containing $\Delta^+_n$. Then it defines a Borel subgroup
$$B=B(\frak{j}, \Delta^+)$$
of $G_\bc$ contained in $Q$.

Let $P$ be a parabolic subgroup of $G_\bc$ containing $B$. Let $S$ be a non-closed $K_\bc$-$P$ double coset in $G_\bc$ and let $S_0$ denote the dense $K_\bc$-$B$ double coset in $S$. By \cite{M1} Theorem 2 we can write
$$S_0=K_\bc c_{\gamma_1}\cdots c_{\gamma_k}wB$$
with some $w\in W$ and a strongly orthogonal system $\{\gamma_1,\ldots,\gamma_k\}$ of roots in $\Delta^+_n$. Here $W$ is the Weyl group of $\Delta$ and
$$c_{\gamma_j}=\exp (X-\overline{X})$$
with some $X\in\frak{g}_\bc(\frak{j},\gamma_j)$ such that $c_{\gamma_j}^2$ is the reflection with respect to $\gamma_j$.

Let $\Theta$ denote the subset of $\Psi$ such that
$$P=BW_\Theta B$$
where $W_\Theta$ is the subgroup of $W$ generated by $\{w_\alpha\mid \alpha\in\Theta\}$. Let $\Delta_\Theta$ denote the subset of $\Delta$ defined by
$$\Delta_\Theta=\{\beta\in\Delta\mid \beta=\sum_{\alpha\in\Theta} n_\alpha \alpha\mbox{ for some }n_\alpha\in\bz\}.$$
If $\gamma_j\in w\Delta_\Theta$ for all $j=1,\ldots,k$, then it follows that
$$c_{\gamma_j}\in wPw^{-1}$$
for all $j=1,\ldots,k$ and therefore
$$Sw^{-1}=S_0Pw^{-1}=K_\bc c_{\gamma_1}\cdots c_{\gamma_k}wPw^{-1}=K_\bc wPw^{-1}$$
becomes closed in $G_\bc$, contradicting the assumption. Hence there exists a $j$ such that $\gamma_j\notin w\Delta_\Theta$. Replacing the order of $\gamma_1,\ldots,\gamma_k$, we may assume that
$$\gamma_1\notin w\Delta_\Theta.$$

Let $\frak{l}$ denote the complex Lie subalgebra of $\frak{g}_\bc$ generated by $\frak{g}_\bc(\frak{j},\gamma_1)\oplus \frak{g}_\bc(\frak{j},-\gamma_1)$ which is isomporphic to $\frak{sl}(2,\bc)$ and let $L$ be the analytic subgroup of $G_\bc$ for $\frak{l}$. Then we have $(L\cap K_\bc)c_{\gamma_1}(L\cap wBw^{-1})=(L\cap K_\bc)c_{\gamma_1}^{-1}(L\cap wBw^{-1})$ since both of the double cosets are open dense in $L$. Hence we have
$$S_0=K_\bc c_{\gamma_1}\cdots c_{\gamma_k}wB= K_\bc c_{\gamma_1}^{-1}c_{\gamma_2}\cdots c_{\gamma_k}wB =K_\bc c_{\gamma_1}\cdots c_{\gamma_k}w_{\gamma_1}wB.$$
If $\gamma_1\notin w\Delta^+$, then $\gamma_1\in w_{\gamma_1}w\Delta^+$. So we may assume
$$\gamma_1\in w\Delta^+$$
replacing $w$ with $w_{\gamma_1}w$ if necessary. Let $\ell$ denote the real rank of $G_\br$.

\begin{lemma} \ There exists a maximal strongly orthgonal system $\{\beta_1,\ldots,\beta_\ell\}$ of roots in $\Delta^+_n$ satisfying the following conditions.

{\rm (i)} \ If $\gamma_1$ is a long root of $\Delta$, then $\beta_1=\gamma_1$ and $\gamma_2,\ldots,\gamma_k\in \br\beta_2\oplus\cdots \oplus\br\beta_\ell$. $($If the roots in $\Delta$ have the same length, then we define that all the roots are long roots.$)$

{\rm (ii)} \ If $\gamma_1$ is a short root of $\Delta$, then $\gamma_1\in \br\beta_1\oplus\br\beta_2$ and $\gamma_2,\ldots,\gamma_k\in \br\beta_3\oplus\cdots \oplus\br\beta_\ell$.
\end{lemma}

\begin{proof} First suppose that $\frak{g}_\br$ is of type AIII, DIII, EIII, EVII or DI(of real rank 2). Then the roots in $\Delta$ have the same length. So we have only to take $\beta_j=\gamma_j$ for $j=1,\ldots,k$ and choose an orthogonal system $\{\beta_1,\ldots,\beta_\ell\}$ of roots in $\Delta^+_n$ containing $\{\beta_1,\ldots,\beta_k\}$.

Next suppose that $\frak{g}_\br\cong \frak{sp}(\ell,\br)$. Write
$$\Delta=\{\pm e_r\pm e_s\mid 1\le r<s\le \ell\}\sqcup \{\pm 2e_r\mid 1\le r\le \ell\}$$
and
$$\Delta^+_n=\{e_r+e_s\mid 1\le r<s\le \ell\}\sqcup \{2e_r\mid 1\le r\le \ell\}$$
as usual using an orthonormal basis $\{e_1,\ldots,e_\ell\}$ of $\frak{j}^*$. If $\gamma_1= 2e_r$, then $\{\beta_2,\ldots,\beta_\ell\}=\{2e_s\mid s\ne r\}$ satisfies the condition (i). If $\gamma_1=e_r+e_s$ with $r\ne s$, then we put $\beta_1=2e_r$ and $\beta_2=2e_s$. The assertion (ii) is clear if we put $\{\beta_3,\ldots,\beta_\ell\}=\{2e_p\mid p\ne r,s\}$.

Finally suppose that $\frak{g}_\br =\frak{so}(2,2p-1)$ with $p\ge 2$. Then the real rank of $\frak{g}_\br$ is two and we can write
$$\Delta=\{\pm e_r\pm e_s\mid 1\le r<s\le p\}\sqcup \{\pm e_r\mid 1\le r\le p\}$$
and
$$\Delta^+_n=\{e_1\pm e_s\mid 2\le s\le p\}\sqcup \{e_1\}$$
with an orthonormal basis $\{e_1,\ldots,e_p\}$ of $\frak{j}^*$. If $k=2$, then we have $\gamma_1=\beta_1=e_1\pm e_s$ and $\gamma_2=\beta_2=e_1\mp e_s$ with some $s$. If $k=1$ and $\gamma_1=e_1\pm e_s$, then $\beta_1=\gamma_1$ and $\beta_2=e_1\mp e_s$. If $k=1$ and $\gamma_1=e_1$, then we may put $\beta_1=e_1+e_2$ and $\beta_2=e_1-e_2$.
\end{proof}

\begin{definition} \ (i) \ Define a subroot system $\Delta_1$ of $\Delta$ as follows.

If $\gamma_1$ is a long root of $\Delta$, then we put
$$\Delta_1=\{\pm\beta_1\}=\{\pm\gamma_1\}.$$
On the other hand if $\gamma_1$ is a short root of $\Delta$, then we put
$$\Delta_1=\Delta\cap (\br\beta_1\oplus\br\beta_2)$$
(which is of type ${\rm C}_2$).

(ii) \ Put $\Delta_2=\{\alpha\in\Delta\mid \alpha\mbox{ is orthogonal to }\Delta_1\}$.

(iii) \ Let $\frak{l}_j$ denote the complex Lie subalgebra of $\frak{g}_\bc$ generated by $\bigoplus_{\alpha\in\Delta_j}\frak{g}_\bc(\frak{j},\alpha)$ for $j=1,2$.

(iv) \ Let $L_1$ and $L_2$ denote the analytic subgroups of $G_\bc$ for $\frak{l}_1$ and $\frak{l}_2$, respectively.
\end{definition}

It follows from Lemma 2.1 that
$$c_{\gamma_1}\in L_1\quad\mbox{and that}\quad c_{\gamma_2}\cdots c_{\gamma_k}\in L_2.$$

Let $X_j$ be nonzero root vectors in $\frak{g}_\bc(\frak{j},\beta_j)$ for $j=1,\ldots,\ell$. Then we can define a maximal abelian subspace
$$\frak{t}=\br(X_1-\overline{X_1})\oplus\cdots\oplus\br(X_\ell-\overline{X_\ell})$$
in $i\frak{m}$ and a maximal abelian subspace
$$\frak{a}=\br(X_1+\overline{X_1})\oplus\cdots\oplus\br(X_\ell+\overline{X_\ell})$$
in $\frak{m}$ as in \cite{GM1} Section 2. Since the restricted root system $\Sigma(\frak{t})$ is of type ${\rm BC}_\ell$, the set $\frak{t}^+$ is defined by the long roots in $\Sigma(\frak{t})$. Hence it is of the form
$$\frak{t}^+=\{Y_1+\cdots+Y_\ell\mid Y_j\in\frak{t}_1^+\}$$
where $\frak{t}_j^+=\{s(X_j-\overline{X_j})\mid -(\pi/4)<s<\pi/4\}$ by a suitable normalization of $X_j$ for $j=1,\ldots,\ell$.

Put $T^+=\exp\frak{t}^+$ and $A=\exp\frak{a}$. Then it is shown in \cite{GM1} Lemma 2.1 that $AQ=T^+Q$ and hence that
$$G_\br Q=KAQ=KT^+Q$$
by the Cartan decomposition $G_\br=KAK$. The closure of $G_\br Q$ in $G_\bc$ is written as
$$(G_\br Q)^{cl}=G_\br Q\sqcup G_\br c_{\beta_1}Q\sqcup G_\br c_{\beta_1}c_{\beta_2}Q\sqcup \cdots\sqcup G_\br c_{\beta_1}\cdots c_{\beta_\ell}Q$$
where $c_{\beta_j}=\exp(\pi/4)(X_j-\overline{X_j})$ for $j=1,\ldots,\ell$ (\cite{WZ0} Theorem 3.8). We can also see that
\begin{equation}
G_\br c_{\beta_1}\cdots c_{\beta_k}Q=Kc_{\beta_1}\cdots c_{\beta_k}T_{k+1}^+\cdots T_\ell^+Q \tag{2.1}
\end{equation}
where $T_j^+=\exp\frak{t}_j^+$ since we can consider the action of the Weyl group $W_K(T)$ on $T$ which is of type ${\rm BC}_\ell$.

By the map
$$\iota: xK_\bc\mapsto (xQ,x\overline{Q})$$
the complex symmetric space $G_\bc/K_\bc$ is embedded in $G_\bc/Q\times G_\bc/\overline{Q}$ (\cite{WZ1}). It is shown in \cite{BHH} Section 3 and \cite{GM1} Proposition 2.2 that
$$\iota(D/K_\bc)=G_\br Q/Q\times G_\br \overline{Q}/\overline{Q}.$$

\begin{lemma} \ Suppose that
$$\iota(xK_\bc)\in G_\br c_{\beta_1}Q/Q\times G_\br\overline{Q}/\overline{Q}$$
and that $\gamma_1$ is a long root of $\Delta_n^+$. $($If the roots in $\Delta$ have the same length, then we define that all the roots are long roots.$)$ Define a $K_\bc$-$B$ double coset $\widetilde{S}_1$ by
$$\widetilde{S}_1=K_\bc c_{\gamma_2}\cdots c_{\gamma_k}wB.$$
Then $\widetilde{S}_1$ is contained in $\p S=S^{cl}-S$ and
$$x\widetilde{S}_1\cap S'_0\ne\phi.$$
\end{lemma}

\begin{proof} It is clear that we may replace $x$ by any elements in the double coset $G_\br xK_\bc$. By the left $G_\br$-action we may assume that $x\in\overline{Q}$. By the right $K_\bc$-action we may moreover assume that $x\in\overline{N}$ since $\overline{Q}=\overline{N}K_\bc$. Since $K=K_\bc\cap G_\br$ normalizes $\overline{N}$, we may assume by (2.1) that
$$xQ=c_{\beta_1}t_2\cdots t_\ell Q$$
with some $t_j\in T^+_j$ for $j=2,\ldots,\ell$. As in \cite{WZ1}, we write
$$c_{\beta_1}=c_{\gamma_1}=c=c^-c^+\mand t_j=t_j^-t_j^+\mbox{ for }j=2,\ldots,\ell$$
with $c^-,\ t_j^-\in\overline{N}$ and $c^+,\ t_j^+\in Q$. Then we have
$$x=c^-t_2^-\cdots t_\ell^-.$$

It follows from Lemma 2.1 and Definition 2.2 that $c_{\gamma_2}\cdots c_{\gamma_k}\in L_2$. Since $\Ad(c_{\gamma_2}\cdots c_{\gamma_k}) \frak{j}$ is $\theta$-stable, the double cosets
$$S_{L_2}=(L_2\cap K_\bc)c_{\gamma_2}\cdots c_{\gamma_k}(L_2\cap wBw^{-1}) \mand S'_{L_2}=(L_2\cap G_\br)c_{\gamma_2}\cdots c_{\gamma_k}(L_2\cap wBw^{-1})$$
correspond by the duality (\cite{M1} Theorem 2). 

It follows from Lemma 2.1 (i) and Definition 2.2 that
$$c^\pm\in L_1\mand t_2^\pm,\ldots, t_\ell^\pm\in L_2.$$
It follows moreover from Definition 2.2 (i) that $\frak{l}_1\cong \frak{sl}(2,\bc)$.

Write $y=t_2^-\cdots t_\ell^-$. Then we have
$$yQ=t_2\cdots t_\ell Q\subset T^+Q\subset G_\br Q$$
and
$$y\overline{Q}=\overline{Q}\subset G_\br\overline{Q}.$$
Hence we have
\begin{equation}
y\in L_2\cap (C(S_1)\cap C(S_2))=L_2\cap D \notag
\end{equation}
by \cite{GM1} (1.3). By the inclusion (1.2) this implies that the set $yS_{L_2}\cap S'_{L_2}$ is nonempty and closed in $L_2$. Take an element $z$ of $yS_{L_2}\cap S'_{L_2}$.

Since $\gamma_1\in w\Delta^+$, we have $c^+\in wBw^{-1}$. Since $c^+\in L_1$ commutes with elements in $L_2$, we have
\begin{align*}
cz\in cyS_{L_2} & = c^-c^+y(L_2\cap K_\bc) c_{\gamma_2}\cdots c_{\gamma_k}(L_2\cap wBw^{-1}) \\
& =c^-y(L_2\cap K_\bc) c_{\gamma_2}\cdots c_{\gamma_k}c^+(L_2\cap wBw^{-1}) \\
& \subset c^-yK_\bc c_{\gamma_2}\cdots c_{\gamma_k}wBw^{-1} =x\widetilde{S}_1w^{-1}
\end{align*}
On the other hand we have
\begin{align*}
cz\in cS'_{L_2} & = c(L_2\cap G_\br) c_{\gamma_2}\cdots c_{\gamma_k}(L_2\cap wBw^{-1}) \\
& =(L_2\cap G_\br) c_{\gamma_1}c_{\gamma_2}\cdots c_{\gamma_k}(L_2\cap wBw^{-1}) \subset S'_0w^{-1}.
\end{align*}
Hence $x\widetilde{S}_1\cap S'_0\ne\phi$. It is clear that $\widetilde{S}_1\subset S_0^{cl}=S^{cl}$ because $(L_1\cap K_\bc)(L_1\cap wBw^{-1})\subset ((L_1\cap K_\bc)c(L_1\cap wBw^{-1}))^{cl}=L_1$.

Now we will prove $\widetilde{S}_1\not\subset S$. Consider the map
$$\varphi:K_\bc\backslash G_\bc/B\ni K_\bc gB\mapsto B\theta(g)^{-1}gB\in B\backslash G_\bc/B$$
introduced in \cite{Sp} where $\theta$ is the holomorphic involution in $G_\bc$ defining $K_\bc$. We have
$$\varphi(\widetilde{S}_1)=Bw^{-1}w_{\gamma_2}\cdots w_{\gamma_k}wB$$
and
$$\varphi(S) =\varphi(S_0P)\subset Pw^{-1}w_{\gamma_1}\cdots w_{\gamma_k}wP =BW_\Theta w^{-1}w_{\gamma_1}\cdots w_{\gamma_k}wW_\Theta B.$$
So we have only to show
\begin{equation}
w^{-1}w_{\gamma_2}\cdots w_{\gamma_k}w\notin W_\Theta w^{-1}w_{\gamma_1}\cdots w_{\gamma_k}wW_\Theta. \tag{2.2}
\end{equation}
Let $Z$ be an element in $\frak{j}$ defining $P$. This implies that $Z$ is dominant for $\Delta^+$ and that $\{\alpha\in\Psi\mid \alpha(Z)=0\}=\Theta$. Let $w_1$ and $w_2$ be elements in $W_\Theta$. Let $B(\ ,\ )$ denote the Killing form on $\frak{g}$ and let $Y_{\gamma_1}$ denote the element in $\frak{j}$ such that
$$\gamma_1(Y)=B(Y,Y_{\gamma_1})\quad\mbox{for all }Y\in\frak{j}.$$
Then we have
\begin{align*}
& \quad B(Z,w^{-1}w_{\gamma_2}\cdots w_{\gamma_k}wZ)- B(Z,w_1w^{-1}w_{\gamma_1}w_{\gamma_2}\cdots w_{\gamma_k}ww_2Z) \\
& =B(wZ-w_{\gamma_1}wZ,w_{\gamma_2}\cdots w_{\gamma_k}wZ) \\
& ={2B(Y_{\gamma_1},wZ)\over B(Y_{\gamma_1},Y_{\gamma_1})}B(Y_{\gamma_1}, w_{\gamma_2}\cdots w_{\gamma_k}wZ) \\
& ={2B(Y_{\gamma_1},wZ)^2\over B(Y_{\gamma_1},Y_{\gamma_1})}>0
\end{align*}
since $\gamma_1\notin w\Delta_\Theta$. Thus we have proved (2.2).
\end{proof}

\begin{lemma} \ Suppose that 
$$\iota(xK_\bc)\in G_\br c_{\beta_1}Q/Q\times G_\br\overline{Q}/\overline{Q}$$
and that $\gamma_1$ is a short root of $\Delta_n^+$. $($We assume that $\frak{g}_\br\cong \frak{sp}(\ell,\br)$ or $\frak{so}(2,2p-1)$ with $p\ge 2$.$)$ Define a $K_\bc$-$B$ double coset $\widetilde{S}_1$ by
$\widetilde{S}_1=K_\bc gc_{\gamma_2}\cdots c_{\gamma_k}wB$
where
$$g=\begin{cases} e & \text{if $\gamma_1$ is the simple short root of $\Delta_1^+$,} \\
c_\beta & \text{if $\gamma_1$ is the non-simple short root of $\Delta_1^+$.} \end{cases}$$
Here $\Delta_1^+=\Delta_1\cap w\Delta^+$ and $\beta$ is the long simple root of $\Delta_1^+$.
Then $\widetilde{S}_1$ is contained in $\p S=S^{cl}-S$ and
$$x\widetilde{S}_1\cap {S'_0}^{cl}\ne\phi.$$
\end{lemma}

\begin{proof}
It follows from Lemma 2.1 (ii) and Definition 2.2 that
$$c_{\beta_1}^\pm,t_2^\pm\in L_1\mand t_3^\pm,\ldots, t_\ell^\pm\in L_2.$$
It follows moreover from Definition 2.2 (i) that $\frak{l}_1\cong \frak{sp}(2,\bc)$.

Write $y=t_3^-\cdots t_\ell^-$. Then by the same argument for long $\gamma_1$ we see that the set $yS_{L_2}\cap S'_{L_2}$ is nonempty and closed in $L_2$. Take an element $z$ of $yS_{L_2}\cap S'_{L_2}$.

The positive system $\Delta_1^+$ of $\Delta_1$ consists of two long roots and two short roots. Since $\gamma_1\in \Delta_1^+$, $\gamma_1$ is either of these two short roots. Write $x_1=c_{\beta_1}^-t_2^-$.

First assume that $\gamma_1$ is the simple short root of $\Delta_1^+$. Then it follows from Proposition 3.2 (i) in the next section that
\begin{equation}
x_1(L_1\cap K_\bc)(L_1\cap wBw^{-1})\cap ((L_1\cap G_\br)c_{\gamma_1}(L_1\cap wBw^{-1}))^{cl} \tag{2.3}
\end{equation}
is nonempty. Note that $L_1\cap wBw^{-1}$ and $\gamma_1$ correspond to $w_{\beta_2}Bw_{\beta_2}^{-1}$ and $\delta$ in the next section, respectively. Let $z_1$ be an element of (2.3). Then we have
\begin{align*}
z_1z &\in x_1(L_1\cap K_\bc)(L_1\cap wBw^{-1})yS_{L_2} \\
&=x_1(L_1\cap K_\bc)(L_1\cap wBw^{-1})y(L_2\cap K_\bc)c_{\gamma_2}\cdots c_{\gamma_k}(L_2\cap wBw^{-1}) \\
&=x_1y(L_1\cap K_\bc)(L_2\cap K_\bc)c_{\gamma_2}\cdots c_{\gamma_k}(L_1\cap wBw^{-1})(L_2\cap wBw^{-1}) \\
&\subset xK_\bc c_{\gamma_2}\cdots c_{\gamma_k}wBw^{-1}=x\widetilde{S}_1w^{-1}
\end{align*}
and
\begin{align*}
z_1z &\in ((L_1\cap G_\br)c_{\gamma_1}(L_1\cap wBw^{-1}))^{cl} S'_{L_2} \\
&= ((L_1\cap G_\br)c_{\gamma_1}(L_1\cap wBw^{-1}))^{cl} (L_2\cap G_\br)c_{\gamma_2}\cdots c_{\gamma_k}(L_2\cap wBw^{-1}) \\
&\subset (G_\br c_{\gamma_1}c_{\gamma_2}\cdots c_{\gamma_k}wBw^{-1})^{cl}={S'_0}^{cl}w^{-1}.
\end{align*}
So we have $x\widetilde{S}_1\cap {S'_0}^{cl}\ne\phi$. We can prove $\widetilde{S}_1\subset S^{cl}-S$ by the same arguments as in the proof of Lemma 2.3.

Next assume that $\gamma_1$ is the non-simple short root of $\Delta_1^+$. Then it follows from Proposition 3.2 (ii) in the next section that
\begin{equation}
x_1(L_1\cap K_\bc)c_\beta(L_1\cap wBw^{-1})\cap ((L_1\cap G_\br)c_{\gamma_1}(L_1\cap wBw^{-1}))^{cl} \notag
\end{equation}
is nonempty. Note that $L_1\cap wBw^{-1},\ \gamma_1$ and $\beta$ correspond to $B,\ \delta$ and $\beta_2$ in the next section, respectively. By the same argument as above we can prove
$$x\widetilde{S}_1\cap {S'_0}^{cl}\ne\phi.$$
It follows from Remark 3.3 that $\widetilde{S}_1\subset S^{cl}$. Finally we will prove that $\widetilde{S}_1\not\subset S$. Using the same argument as in the proof of Lemma 2.3, we have only to show
\begin{equation}
w^{-1}w_\beta w_{\gamma_2}\cdots w_{\gamma_k}w\notin W_\Theta w^{-1}w_{\gamma_1}\cdots w_{\gamma_k}wW_\Theta. \tag{2.4}
\end{equation}
Let $Z$ and $Y_{\gamma_1}$ be as  in the proof of Lemma 2.3. Define $Y_\beta\in \frak{j}$ so that
$$\beta(Y)=B(Y,Y_\beta)\quad\mbox{for all }Y\in\frak{j}.$$
Then we have
\begin{align*}
& \quad B(Z,w^{-1}w_\beta w_{\gamma_2}\cdots w_{\gamma_k}wZ)- B(Z,w_1w^{-1}w_{\gamma_1}w_{\gamma_2}\cdots w_{\gamma_k}ww_2Z) \\
& =B(w_\beta wZ-w_{\gamma_1}wZ,w_{\gamma_2}\cdots w_{\gamma_k}wZ) \\
& =B(wZ-w_{\gamma_1}wZ,w_{\gamma_2}\cdots w_{\gamma_k}wZ)
-B(wZ-w_\beta wZ,w_{\gamma_2}\cdots w_{\gamma_k}wZ) \\
& ={2B(Y_{\gamma_1},wZ)\over B(Y_{\gamma_1},Y_{\gamma_1})}B(Y_{\gamma_1}, w_{\gamma_2}\cdots w_{\gamma_k}wZ) -{2B(Y_\beta,wZ)\over B(Y_\beta,Y_\beta)}B(Y_\beta, w_{\gamma_2}\cdots w_{\gamma_k}wZ)  \\
& ={2B(Y_{\gamma_1},wZ)^2\over B(Y_{\gamma_1},Y_{\gamma_1})} -{2B(Y_\beta,wZ)^2\over B(Y_\beta,Y_\beta)} >0
\end{align*}
for $w_1,w_2\in W_\Theta$ since
$$B(Y_{\gamma_1},wZ)>0,\quad 0\le B(Y_\beta,wZ)\le B(Y_{\gamma_1},wZ)\mand B(Y_\beta,Y_\beta)=2B(Y_{\gamma_1},Y_{\gamma_1}).$$
Thus we have proved (2.4).
\end{proof}

Using the conjugation on $G_\bc$ with respect to the real form $G_\br$, it follows from Lemma 2.3 and Lemma 2.4 the following.

\begin{corollary} \ Suppose that 
$$\iota(xK_\bc)\in G_\br Q/Q\times G_\br\overline{c_{\beta_1}}\overline{Q}/\overline{Q}.$$
Then there exists a $K_\bc$-$B$ double coset $\widetilde{S}_2$ contained in $\p S$ such that
$$x\widetilde{S}_2\cap {S'_0}^{cl}\ne\phi.$$
\end{corollary}

\noindent {\it Proof of Theorem 1.3}. \ Let $S$ be a non-closed $K_\bc$-$P$ double coset in $G_\bc$. Then it follows from Lemma 2.3, Lemma 2.4 and Corollary 2.5 that there exist $K_\bc$-$B$ double cosets $\widetilde{S}_1$ and $\widetilde{S}_2$ contained in $\p S$ such that
\begin{equation}
x(\widetilde{S}_1\cup \widetilde{S}_2)\cap {S'_0}^{cl}\ne\phi \tag{2.5}
\end{equation}
for all $x\in\p D$ satisfying
\begin{equation}
xK_\bc \in \iota^{-1}((G_\br c_{\beta_1}Q/Q\times G_\br\overline{Q}/\overline{Q}) \sqcup (G_\br Q/Q\times G_\br\overline{c_{\beta_1}}\overline{Q}/\overline{Q})). \tag{2.6}
\end{equation}
Suppose that
$$y(\widetilde{S}_1\cup \widetilde{S}_2)^{cl}\cap {S'_0}^{cl}=\phi.$$
for some $y\in\p D$. Then there exists a neighborhood $U$ of $y$ in $G_\bc$ such that
$$x(\widetilde{S}_1\cup \widetilde{S}_2)^{cl}\cap {S'_0}^{cl}=\phi$$
for all $x\in U$. But this contradicts (2.5) because the right hand side of (2.6) is dense in $\p(D/K_\bc)$. \hfill $\square$

\section{$Sp(2,\br)$-case}

Let $G_\bc=Sp(2,\bc)=\{g\in GL(4,\bc)\mid {}^tgJg=J\}$ where
$$J=\bp 0 & -I_2 \\ I_2 & 0 \ep.$$
Let
$$K_\bc=\left\{\bp g & 0 \\ 0 & {}^tg^{-1} \ep\Bigm| g\in GL(2,\bc)\right\}
\mand
G_\br=G_\bc\cap U(2,2)\cong Sp(2,\br).$$
Put $U_+=\bc e_1\oplus\bc e_2$ and $U_-=\bc e_3\oplus\bc e_4$ by using the canonical basis $\{e_1,e_2,e_3,e_4\}$ of $\bc^4$. Then
$$K_\bc=Q\cap\overline{Q}$$
where $Q=\{g\in G_\bc\mid gU_+=U_+\}$ and $\overline{Q}=\{g\in G_\bc\mid gU_-=U_-\}$. (Here $\overline{*}$ is the conjugate of $*$ with respect to the real form $G_\br$ of $G_\bc$.)

The full flag manifold $X$ of $G_\bc$ consists of the flags
$$(V_1,V_2)$$
in $\bc^4$ where $\dim V_j=j, \ V_1\subset V_2$ and ${}^tuJv=0$ for all $u,v\in V_2$. Let $B$ denote the Borel subgroup of $G_\bc$ defined by
$$B=\{g\in G_\bc\mid g\bc e_1=\bc e_1\mbox{ and }gU_+=U_+\}.$$
Then the full flag manifold $X$ is identified with $G_\bc/B$ by the map
$$gB\mapsto (V_1,V_2)=(g\bc e_1,gU_+).$$

There are eleven $K_\bc$-orbits
\begin{align*}
S_1 &=\{(V_1,V_2)\mid V_2=U_+\}, \\
S_2 &=\{(V_1,V_2)\mid V_2=U_-\}, \\
S_3 &=\{(V_1,V_2)\mid V_1\subset U_+,\ \dim(V_2\cap U_-)=1\}, \\
S_4 &=\{(V_1,V_2)\mid V_1\subset U_-,\ \dim(V_2\cap U_+)=1\}, \\
S_5 &=\{(V_1,V_2)\mid V_1\subset U_+\}-(S_1\sqcup S_3), \\
S_6 &=\{(V_1,V_2)\mid V_1\subset U_-\}-(S_2\sqcup S_4), \\
S_7 &=\{(V_1,V_2)\mid \dim(V_2\cap U_+)=\dim(V_2\cap U_-)=1\}-(S_3\sqcup S_4), \\
S_8 &=\{(V_1,V_2)\mid V_1\cap U_+=\{0\},\ \dim(V_2\cap U_+)=1,\ V_2\cap U_-=\{0\}\}, \\
S_9 &=\{(V_1,V_2)\mid V_1\cap U_-=\{0\},\ \dim(V_2\cap U_-)=1,\ V_2\cap U_+=\{0\}\}, \\
S_{10} &=\{(V_1,V_2)\mid V_2\cap U_\pm =\{0\},\ {}^tvJ\tau(v)=0\mbox{ for }v\in V_1\}, \\
S_\op &=\{(V_1,V_2)\mid V_2\cap U_\pm =\{0\},\ {}^tvJ\tau(v)\ne 0\mbox{ for }v\in V_1-\{0\}\}
\end{align*}
on $X$ where
$$\tau(v)=\bp I_2 & 0 \\ 0 & -I_2 \ep v$$
for $v\in\bc^4$. These orbits are related as follows (\cite{MO} Fig. 12).

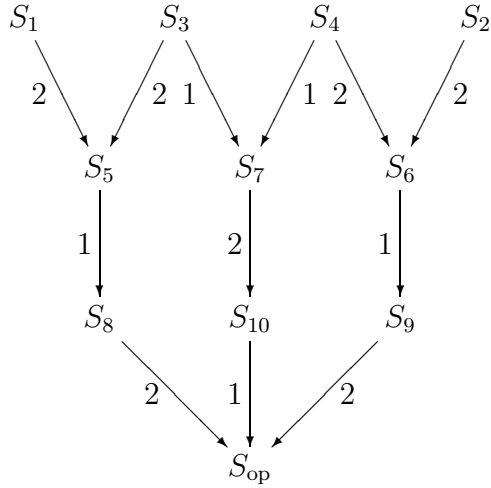
\setlength{\unitlength}{1mm}
\begin{picture}(100,70)
\put(40,65){\makebox(0,0){$S_1$}}
\put(60,65){\makebox(0,0){$S_3$}}
\put(80,65){\makebox(0,0){$S_4$}}
\put(100,65){\makebox(0,0){$S_2$}}
\put(50,45){\makebox(0,0){$S_5$}}
\put(70,45){\makebox(0,0){$S_7$}}
\put(90,45){\makebox(0,0){$S_6$}}
\put(50,25){\makebox(0,0){$S_8$}}
\put(70,25){\makebox(0,0){$S_{10}$}}
\put(90,25){\makebox(0,0){$S_9$}}
\put(70,5){\makebox(0,0){$S_\op$}}
\put(41.5,62){\vector(1,-2){7}}
\put(58.5,62){\vector(-1,-2){7}}
\put(61.5,62){\vector(1,-2){7}}
\put(78.5,62){\vector(-1,-2){7}}
\put(81.5,62){\vector(1,-2){7}}
\put(98.5,62){\vector(-1,-2){7}}
\put(50,42){\vector(0,-1){14}}
\put(70,42){\vector(0,-1){14}}
\put(90,42){\vector(0,-1){14}}
\put(53,22){\vector(1,-1){14}}
\put(87,22){\vector(-1,-1){14}}
\put(70,22){\vector(0,-1){14}}
\put(42,55){\makebox(0,0){2}}
\put(58,55){\makebox(0,0){2}}
\put(62,55){\makebox(0,0){1}}
\put(78,55){\makebox(0,0){1}}
\put(82,55){\makebox(0,0){2}}
\put(98,55){\makebox(0,0){2}}
\put(48,35){\makebox(0,0){1}}
\put(68,35){\makebox(0,0){2}}
\put(88,35){\makebox(0,0){1}}
\put(57,15){\makebox(0,0){2}}
\put(83,15){\makebox(0,0){2}}
\put(68,15){\makebox(0,0){1}}
\end{picture}

Let $P_1$ and $P_2$ be the parabolic subgroups of $G_\bc$ defined by
$$P_1=Q\mand P_2=\{g\in G_\bc\mid g\bc e_1=\bc e_1\},$$
respectively. Then the above diagram implies, for example, that
$$S_1P_2=S_5P_2\quad\mbox{and that}\quad \dim S_1=\dim S_5-1$$
by the arrow attached with the number 2 joining $S_1$ and $S_5$.

On the other hand define subsets
$$C_+=\{z\in\bc^4\mid (z,z)>0\},\quad C_-=\{z\in\bc^4\mid (z,z)<0\}$$
$$\mand C_0=\{z\in\bc^4\mid (z,z)=0\}$$
of $\bc^4$ using the Hermitian form $(w,z)= \overline{w_1}z_1+\overline{w_2}z_2 -\overline{w_3}z_3 -\overline{w_4}z_4$ defining $U(2,2)$. For $v\in \bc^4$ define subspaces
$$v^J=\{u\in \bc^4\mid {}^tvJu=0\}\mand v^\perp=\{u\in \bc^4\mid (v,u)=0\}$$
of $\bc^4$. Then $C_0$ is devided as $C_0=C_0^s\sqcup C_0^r$ where
$$C_0^s=\{v\in C_0\mid v^J=v^\perp\}\mand C_0^r=\{v\in C_0\mid v^J\ne v^\perp\}.$$

The $G_\br$-orbits on $X$ are
\begin{align*}
S'_1 &=\{(V_1,V_2)\mid V_2-\{0\}\subset C_+\}, \\
S'_2 &=\{(V_1,V_2)\mid V_2-\{0\}\subset C_-\}, \\
S'_3 &=\{(V_1,V_2)\mid V_1-\{0\}\subset C_+,\ V_2\cap C_-\ne\phi\}, \\
S'_4 &=\{(V_1,V_2)\mid V_1-\{0\}\subset C_-,\ V_2\cap C_+\ne\phi\}, \\
S'_5 &=\{(V_1,V_2)\mid V_1-\{0\}\subset C_+,\ V_2\cap C_0^s\ne\{0\}\}, \\
S'_6 &=\{(V_1,V_2)\mid V_1-\{0\}\subset C_-,\ V_2\cap C_0^s\ne\{0\}\}, \\
S'_7 &=\{(V_1,V_2)\mid V_1-\{0\}\subset C_0^r,\ V_2\not\subset C_0\}, \\
S'_8 &=\{(V_1,V_2)\mid V_1\subset C_0^s,\ V_2\cap C_+\ne\phi\}, \\
S'_9 &=\{(V_1,V_2)\mid V_1\subset C_0^s,\ V_2\cap C_-\ne\phi\}, \\
S'_{10} &=\{(V_1,V_2)\mid V_1-\{0\}\subset C_0^r,\ V_2\subset C_0\}, \\
S'_\op &=\{(V_1,V_2)\mid V_1\subset C_0^s,\ V_2\subset C_0\}.
\end{align*}
Here the $K_\bc$-orbit $S_j$ and the $G_\br$-orbit $S'_j$ correspond by the duality for each $j=1,\ldots,10,\op$.

Take a maximal abelian subspace
$$\frak{j}=\left\{Y(a_1,a_2)=\bp a_1 & 0 & 0 & 0 \\
0 & a_2 & 0 & 0 \\
0 & 0 & -a_1 & 0 \\
0 & 0 & 0 & -a_2 \ep
\Bigm| a_1,a_2\in\br\right\}$$
of $i\frak{m}$. Using the linear forms $e_j:Y(a_1,a_2)\mapsto a_j$ for $j=1,2$, we can write
$$\Delta=\{\pm 2e_1,\pm 2e_2,\pm e_1\pm e_2\}
\mand
\Delta_n^+=\{2e_1,2e_2,e_1+e_2\}.$$
Write $\beta_1=2e_1,\ \beta_2=2e_2$ and $\delta=e_1+e_2$. Take root vectors $X_1=-E_{13}$ of $\frak{g}_\bc(\frak{j},\beta_1)$ and $X_2=-E_{24}$ of $\frak{g}_\bc(\frak{j},\beta_2)$ where $E_{ij}\ (i,j=1,\ldots,4)$ denote the matrix units. Define
$$t_1(s)=\exp s(X_1-\overline{X_1})=\exp s(E_{31}-E_{13}) =\bp
\cos s & 0 & -\sin s & 0 \\
0 & 1 & 0 & 0 \\
\sin s & 0 & \cos s & 0 \\
0 & 0 & 0 & 1 \ep$$
and
$$t_2(s)=\exp s(X_2-\overline{X_2})=\exp s(E_{42}-E_{24})=\bp
1 & 0 & 0 & 0 \\
0 & \cos s & 0 & -\sin s \\
0 & 0 & 1 & 0 \\
0 & \sin s & 0 & \cos s \ep$$
for $s\in\br$. Then we can write the Akhiezer-Gindikin domain $D$ as
$$D=G_\br T^+K_\bc$$
where $T^+=\{t_1(s_1)t_2(s_2)\mid |s_1|<\pi/4,\ |s_2|<\pi/4\}$. Write $c_{\beta_j}=t_j(\pi/4)$ and $w_{\beta_j}=t_j(\pi/2)$ for $j=1,2$. Then we can write
$$S_j=K_\bc gB\mand S'_j=G_\br gB$$
for $j=1,\ldots,10,\op$ with the following representatives $g$ (\cite{M1} Theorem 2).

\bigskip
\centerline{
\vbox{\offinterlineskip
\hrule
\halign{&\vrule#&\strut\ $\hfil#\hfil$\ \cr
height2pt&\omit&&\omit&&\omit&&\omit&&\omit&&\omit&&\omit&&\omit&&\omit&&\omit&&\omit&&\omit& \cr
& j && 1 && 2 && 3 && 4 && 5 && 6 && 7 && 8 && 9 && 10 && \op &\cr
height2pt&\omit&&\omit&&\omit&&\omit&&\omit&&\omit&&\omit&&\omit&&\omit&&\omit& &\omit&&\omit&\cr
\noalign{\hrule}
height4pt&\omit&&\omit&&\omit&&\omit&&\omit&&\omit&&\omit&&\omit&&\omit&&\omit& &\omit&&\omit&\cr
& g && e && w_{\beta_1}w_{\beta_2} && w_{\beta_2} && w_{\beta_1} && c_{\beta_2} && c_{\beta_2}w_{\beta_1} && c_\delta w_{\beta_2} && c_{\beta_1} && c_{\beta_1}w_{\beta_2} && c_\delta && c_{\beta_1}c_{\beta_2} &\cr
height4pt&\omit&&\omit&&\omit&&\omit&&\omit&&\omit&&\omit&&\omit&&\omit&&\omit& &\omit&&\omit&\cr}
\hrule}
}

\noindent Here
$$c_\delta={1\over\sqrt{2}}\bp 1 & 0 & 0 & -1 \\
0 & 1 & -1 & 0 \\
0 & 1 & 1 & 0 \\
1 & 0 & 0 & 1 \ep =\exp {\pi\over 4}(X_\delta-\overline{X_\delta})$$
with $X_\delta=-(E_{14}+E_{23})\in\frak{g}_\bc(\frak{j},\delta)$.

The standard maximal flag manifold $G_\bc/Q$ is identified with the space $Y$ of two dimensional subspaces $V_+$ of $\bc^4$ such that ${}^tuJv=0$ for all $u,v\in V_+$ by the map
$$G_\bc/Q\ni gQ\mapsto V_+=gU_+\in Y.$$
Similarly we also identify $G_\bc/\overline{Q}$ with $Y$ by the map
$$G_\bc/\overline{Q}\ni g\overline{Q}\mapsto V_-=gU_-\in Y.$$

As in Section 2 the complex symmetric space $G_\bc/K_\bc$ is naturally identified with the open subset
$$\{(V_+,V_-)\in G_\bc/Q\times G_\bc/\overline{Q} \mid V_+\cap V_-=\{0\}\}$$
of $G_\bc/Q\times G_\bc/\overline{Q}\cong Y\times Y$ by the map
$$\iota: gK_\bc\mapsto (V_+,V_-)=(gU_+,gU_-).$$
Then the Akhiezer-Gindikin domain $D/K_\bc$ is identified with
$$G_\br Q/Q\times G_\br \overline{Q} /\overline{Q}=\{(V_+,V_-)\in Y\times Y\mid V_+-\{0\}\subset C_+\mbox{ and }V_--\{0\}\subset C_-\}.$$

Let $xK_\bc$ be an element of $\p (D/K_\bc)$ such that $\iota(xK_\bc)\in G_\br c_{\beta_1}Q/Q\times G_\br \overline{Q} /\overline{Q}$. Then it follows from Lemma 2.3 that
$$xK_\bc gB\cap G_\br c_{\beta_1}gB\ne\phi$$
for $g=e,w_{\beta_2}$ and $c_{\beta_2}$. This implies that
\begin{equation}
xS_1\cap S'_8\ne\phi, \tag{3.1}
\end{equation}
\begin{equation}
xS_3\cap S'_9\ne\phi \tag{3.2}
\end{equation}
and that
\begin{equation}
xS_5\cap S'_\op\ne\phi. \tag{3.3}
\end{equation}
Since ${S'_7}^{cl}=\{(V_1,V_2)\mid V_1\subset C_0\}\supset S'_9$, it follows from (3.2) that
\begin{equation}
xS_3\cap {S'_7}^{cl}\ne\phi. \tag{3.4}
\end{equation}
On the other hand since ${S'_{10}}^{cl}\supset S'_\op$, it follows from (3.3) that
\begin{equation}
xS_5\cap {S'_{10}}^{cl}\ne\phi. \tag{3.5}
\end{equation}

\begin{remark} \ {\rm (i)} \ If $\iota(xK_\bc)\in G_\br Q/Q\times G_\br \overline{c_{\beta_1}}\overline{Q}/\overline{Q}$, then we can prove
$$xS_2\cap S'_9\ne\phi,\quad xS_4\cap S'_8\ne\phi,\quad xS_6\cap S'_\op\ne\phi,$$
$$xS_4\cap {S'_7}^{cl}\ne\phi\mand xS_6\cap {S'_{10}}^{cl}\ne\phi$$
in the same way.

{\rm(ii)} \ If we apply \cite{M8} Theorem 1.3 to this case, then we have
$$x\in\p D\Longrightarrow x(S_5\sqcup S_6)^{cl}\cap S'_\op\ne\phi.$$
So we see that the results in this paper are a refinement of this theorem for Hermitian cases.
\end{remark}

By (3.4) and (3.5) we proved the following.

\begin{proposition} \ If $\iota(xK_\bc)\in G_\br c_{\beta_1}Q/Q\times G_\br \overline{Q} /\overline{Q}$. Then we have$:$

{\rm (i)} \ $xK_\bc w_{\beta_2}B\cap (G_\br c_\delta w_{\beta_2}B)^{cl}\ne\phi$.

{\rm (ii)} \ $xK_\bc c_{\beta_2}B\cap (G_\br c_\delta B)^{cl}\ne\phi$.
\end{proposition}

\begin{remark} \ It is clear that $K_\bc w_{\beta_2}B=S_3\subset S_7^{cl}= (K_\bc c_\delta w_{\beta_2}B)^{cl}$ and that $K_\bc c_{\beta_2}B=S_5\subset S_{10}^{cl}= (K_\bc c_\delta B)^{cl}$.
\end{remark}

\end{document}